\documentclass[a4paper,english]{article}
\usepackage[T1]{fontenc}
\usepackage[utf8]{luainputenc}
\usepackage{amsthm}
\usepackage{amsmath}
\usepackage{amssymb}
\usepackage{esint}

\makeatletter

\pdfpageheight\paperheight
\pdfpagewidth\paperwidth

\theoremstyle{plain}
\newtheorem{thm}{\protect\theoremname}[section]
\theoremstyle{definition}
\newtheorem{defn}[thm]{\protect\definitionname}
\theoremstyle{plain}
\newtheorem{lem}[thm]{\protect\lemmaname}
\ifx\proof\undefined
\newenvironment{proof}[1][\protect\proofname]{\par
\normalfont\topsep6\p@\@plus6\p@\relax
\trivlist
\itemindent\parindent
\item[\hskip\labelsep
\scshape
#1]\ignorespaces
}{%
\endtrivlist\@endpefalse
}
\providecommand{\proofname}{Proof}
\fi

\numberwithin{equation}{section}
\date{}

\makeatother

\usepackage{babel}
\providecommand{\definitionname}{Definition}
\providecommand{\lemmaname}{Lemma}
\providecommand{\theoremname}{Theorem}

\begin{document}

\title{Reconstruction of penetrable obstacles in the anisotropic acoustic
scattering}

\author{Yi-Hsuan Lin}
\maketitle
\begin{abstract}
We develop reconstruction schemes to determine penetrable obstacles
in a region of $\mathbb{R}^{2}$ or $\mathbb{R}^{3}$ and we consider
anisotropic elliptic equations. This algorithm uses oscillating-decaying
solutions to the equation. We apply the oscillating-decaying solutions
and the Runge approximation property to the inverse problem of identifying
an inclusion in an anisotropic elliptic differential equation.

\end{abstract}

\section{Introduction}

The special type solutions for elliptic equations or systems play
an essential role in inverse problems since the pioneer work of Cald$\acute{e}$ron.
In \cite{Sylvester}, Sylvester and Uhlmann used complex geometric
optics (CGO) solutions to solve the inverse boundary value problems
for the conductivity equation. Based on CGO solutions, Ikehata proposed
the so called enclosure method to reconstruct the inclusion obstacle,
see \cite{Ike-enclosue}.There are many results in this reconstruction
algorithm, in \cite{UhlWang-disCGO}, they construct CGO-solutions
with polynomial-type phase function for the Helmholtz equation $\Delta u+k^{2}u=0$
or elliptic system having the Laplacian as the principal part. In
\cite{Nakamura}, he constructed a very special solution of a conductivity
equation $\nabla\cdot(\gamma(x)\nabla u)=0$ (called the oscillating-decaying
solutions), the leading parts is also isotropic. However, when the
medium is anisotropic, we need to consider more general elliptic equations,
such as anisotropic scalar elliptic equations $\nabla\cdot(A^{0}(x)\nabla u)+k^{2}u=0$,
where $A^{0}(x)=(a_{ij}^{0}(x))$, $a_{ij}^{0}(x)=a_{ji}^{0}(x)$
and assume the uniform ellipticity condition, that is, $\forall\xi=(\xi_{1},\xi_{2},\cdots\xi_{n})\in\mathbb{R}^{n}$,
$\lambda^{0}|\xi|^{2}\leq\sum_{i,j}a_{ij}^{0}(x)\xi_{i}\xi_{j}\leq\Lambda^{0}|\xi|^{2}$.
In this paper, we want to use the oscillating-decaying solutions in
our reconstruction algorithm. We have some assumptions. First, we
consider this problem in $\mathbb{R}^{3}$ and assume that $D$ is
an unknown obstacle such that $D\Subset\Omega\subset\mathbb{R}^{3}$
with an inhomogeneous index of refraction subset of a larger domain
$\Omega$.and $D$, $\Omega$ are $C^{1}$ domains. Second, we assume
$a_{ij}(x)=a_{ij}^{0}(x)\chi_{\Omega\backslash D}+\widetilde{a_{ij}}(x)\chi_{D}$,
where $\widetilde{a_{ij}}(x)$ is regarded as a perturbation in the
unknown obstacle $D$ and $\widetilde{a_{ij}}(x)$ satisfies $\widetilde{\lambda}|\xi|^{2}\leq\sum_{i,j}\widetilde{a_{ij}}(x)\xi_{i}\xi_{j}\leq\widetilde{\Lambda}|\xi|^{2}$.
Moreover, we need to assume that there exists a universal constant
$0<\widehat{\lambda}\leq\widehat{\Lambda}$ such that $\forall\xi\in\mathbb{R}^{3}$,
we have $\widehat{\lambda}|\xi|^{2}\leq\sum(\widetilde{a_{ij}}(x)\chi_{D}-a_{ij}^{0}(x))\xi_{i}\xi_{j}\leq\widehat{\Lambda}|\xi|^{2}$,
which mean the perturbed term $\widetilde{A}(x)$ is ``greater''
than the unperturbed term $A^{0}$ inside the unknown obstacle $D$.
Denote $A(x)=(a_{ij}(x))$, $A^{0}(x)=(a_{ij}^{0}(x))$ and let $k>0$
and consider the steady state anisotropic acoustic wave equation in
with Dirichlet boundary condition 
\begin{equation}
\begin{cases}
\nabla\cdot(A(x)\nabla u)+k^{2}u=0 & \mbox{ in }\Omega\\
u=f & \mbox{ on }\partial\Omega.
\end{cases}\label{eq:1.1}
\end{equation}
In the unperturbed case, we have 
\begin{equation}
\begin{cases}
\nabla\cdot(A^{0}(x)\nabla u_{0})+k^{2}u_{0}=0 & \mbox{ in }\Omega\\
u_{0}=f & \mbox{ on }\partial\Omega.
\end{cases}\label{eq:1.2}
\end{equation}
In this paper, we assume that $k^{2}$ is not a Dirichlet eigenvalue
of the operator $-\nabla\cdot(A\nabla\bullet)$ and $-\nabla\cdot(A^{0}\nabla\bullet)$
in $\Omega$. It is known that for any $f\in H^{1/2}(\partial\Omega)$,
there exists a unique solution $u$ to (\ref{eq:1.1}). We define
the Dirichlet-to-Neumann map in the anisotropic case, say $\Lambda_{D}:H^{1/2}(\partial\Omega)\to H^{-1/2}(\partial\Omega)$
as the following.
\begin{defn}
$\Lambda_{D}f:=A\nabla u\cdot\nu=\sum_{i.j=1}^{3}a_{ij}\partial_{j}u\cdot\nu_{i}$
and $\Lambda_{\emptyset}f:=A^{0}\nabla u_{0}\cdot\nu=\sum_{i.j=1}^{3}a_{ij}\partial_{j}u_{0}\cdot\nu_{i}$,
where $\nu=(\nu_{1},\nu_{2},\nu_{3})$ is an outer normal on $\partial\Omega$.
\end{defn}
\textbf{Inverse problem}: Identify the location and the convex hull
of $D$ from the DN-map $\Lambda_{D}$. The domain $D$ can also be
treated as an inclusion embedded in $\Omega$. The aim of this work
is to give a reconstruction algorithm for this problem. Note that
the information on the medium parameter $(\widetilde{a_{ij}}(x))$
inside $D$ is not known a priori.

The main tool in our reconstruction method is the oscillating-decaying
solutions for the second order anisotropic elliptic differential equations.
We use the results coming from the paper \cite{JN-ods} to construct
the oscillating-decaying solution. In section 2, we will construct
the oscillating-decaying solutions for anisotropic elliptic equations,
note that even if $k=0$, which means the equation is $\nabla\cdot(A(x)\nabla u)=0$,
we do not have any CGO-type solutions. Roughly speaking, given a hyperplane,
an oscillating-decaying solution is oscillating very rapidly along
this plane and decaying exponentially in the direction transversely
to the same plane. They are also CGO-solutions but with the imaginary
part of the phase function non-negative. Note that the domain of the
oscillating-decaying solutions is not over the whole $\Omega$, so
we need to extend such solutions to the whole domain. Fortunately,
the Runge approximation property provides us a good approach to extend
this special solution in section 3. 

In Ikehata's work, the CGO-solutions are used to define the indicator
function (see \cite{Ike-enclos1} for the definition). In order to
use the oscillating-decaying solutions to the inverse problem of identifying
an inclusion, we have to modify the definition of the indicator function
using the Runge approximation property. It was first recognized by
Lax \cite{Lax} that the Runge approximation property is a consequence
of the weak unique continuation property. In our case, it is clear
that the anisotropic elliptic equation has the weak unique continuation
property if the leading part is Lipschitz continuous.

\section{Construction of oscillating-decaying solutions}

In this section, we follow the paper \cite{JN-ods} to construct the
oscillating-decaying solution in the anisotropic elliptic equations.
In our case, since we only consider a scalar elliptic equation, it's
construction is simpler than the construction in \cite{JN-ods}. Consider
the Dirichlet problem
\begin{equation}
\begin{cases}
\nabla\cdot(A(x)\nabla u)+k^{2}u=0 & \mbox{ in }\Omega\\
u=f & \mbox{ on }\partial\Omega.
\end{cases}\label{eq:2.1}
\end{equation}
 Note that the oscillating-decaying solutions of 
\[
\begin{cases}
\nabla\cdot(A(x)\nabla u)=0 & \mbox{ in }\Omega\\
u=f & \mbox{ on }\partial\Omega
\end{cases}
\]
will have the same representation as the equation (\ref{eq:2.1}),
that is, the lower order term $k^{2}u$ will not affect the form of
the oscillating-decaying solutions, we will see the detail in the
following constructions. Now, we assume that the domain $\Omega$
is an open, bounded smooth domain in $\mathbb{R}^{3}$ and the coefficients
$A(x)=(a_{ij}(x))$ satisfying $\sum_{i.j=1}^{3}a_{ij}(x)\xi_{i}\xi_{j}\geq\lambda|\xi|^{2}$,
$\forall\xi=(\xi_{1},\xi_{2},\xi_{3})\in\mathbb{R}^{3}$ and $\lambda$
is a universal constant.

Assume that 
\[
A(x)=(a_{ij}(x))\in B^{\infty}(\mathbb{R}^{3})=\{f\in C^{\infty}(\mathbb{R}^{3}):\partial^{\alpha}f\in L^{\infty}(\mathbb{R}^{3}),\mbox{ }\forall\alpha\in\mathbb{Z}_{+}^{3}\}
\]
 is the anisotropic coefficients satisfying $a_{ij}(x)=a_{ji}(x)$
$\forall i,j$ and there exists a $\lambda>0$ such that $\sum_{i,j}a_{ij}(x)\xi_{i}\xi_{j}\geq\lambda|\xi|^{2}$
$\forall x\in\mathbb{R}^{3}$ (uniform ellipticity). It is clear that
$A(x)$ is Lipschitz continuous if each $a_{ij}(x)\in B^{\infty}(\mathbb{R}^{3})$,
it has weak continuation property.

We give several notations as follows. Assume that $\Omega\subset\mathbb{R}^{3}$
is an open set with smooth boundary and $\omega\in S^{2}$ is given.
Let $\eta\in S^{2}$ and $\zeta\in S^{2}$ be chosen so that $\{\eta,\zeta,\omega\}$
forms an orthonormal system of $\mathbb{R}^{3}$. We then denote $x'=(x\cdot\eta,x\cdot\zeta)$.
Let $t\in\mathbb{R}$, $\Omega_{t}(\omega)=\Omega\cap\{x\cdot\omega>t\}$
and $\Sigma_{t}(\omega)=\Omega\cap\{x\cdot\omega=t\}$ be a non-empty
open set. We consider a scalar function $u_{\chi_{t},t,b,N,\omega}(x,\tau):=u(x,\tau)\in C^{\infty}(\overline{\Omega_{t}(\omega)}\backslash\partial\Sigma_{t}(\omega))\cap C^{0}(\overline{\Omega_{t}(\omega)})$
with $\tau\gg1$ satisfying:
\begin{equation}
\begin{cases}
L_{A}u=\nabla\cdot(A(x)\nabla u)+k^{2}u=0 & \mbox{ in }\Omega_{t}(\omega)\\
u=e^{i\tau x\cdot\xi}\{\chi_{t}(x')Q_{t}(x')b+\beta_{\chi_{t},t,b,N,\omega}\} & \mbox{ on }\Sigma_{t}(\omega),
\end{cases}\label{eq:2.2}
\end{equation}
where $\xi\in S^{2}$ lying in the span of $\eta$ and $\zeta$ is
chosen and fixed, $\chi_{t}(x')\in C_{0}^{\infty}(\mathbb{R}^{2})$
with supp$(\chi_{t})\subset\Sigma_{t}(\omega)$, $Q_{t}(x')$ is a
nonzero smooth function and $0\neq b\in\mathbb{C}$. Moreover, $\beta_{\chi_{t},b,t,N,\omega}(x',\tau)$
is a smooth function supported in supp($\chi_{t}$) satisfying:
\[
\|\beta_{\chi_{t},b,t,N,\omega}(\cdot,\tau)\|_{L^{2}(\mathbb{R}^{2})}\leq c\tau^{-1}
\]
for some constant $c>0$. From now on, we use $c$ to denote a general
positive constant whose value may vary from line to line. As in the
paper \cite{JN-ods}, $u_{\chi_{t},b,t,N,\omega}$ can be written
as 
\[
u_{\chi_{t},b,t,N,\omega}=w_{\chi_{t},b,t,N,\omega}+r_{\chi_{t},b,t,N,\omega}
\]
 with 
\begin{equation}
w_{\chi_{t},b,t,N,\omega}=\chi_{t}(x')Q_{t}e^{i\tau x\cdot\xi}e^{-\tau(x\cdot\omega-t)A_{t}(x')}b+\gamma_{\chi_{t},b,t,N,\omega}(x,\tau)\label{eq:2.3}
\end{equation}
and $r_{\chi_{t}b,t,N,\omega}$ satisfying 
\begin{equation}
\|r_{\chi_{t},b,t,N,\omega}\|_{H^{1}(\Omega_{t}(\omega))}\leq c\tau^{-N-1/2},\label{eq:2.4}
\end{equation}
where $A_{t}(\cdot)\in B^{\infty}(\mathbb{R}^{2})$ is a complex function
with its real part Re$A_{t}(x')>0$, and $\gamma_{\chi_{t},b,t,N,\omega}$
is a smooth function supported in supp($\chi_{t}$) satisfying 
\begin{equation}
\|\partial_{x}^{\alpha}\gamma_{\chi_{t},b,t,N,\omega}\|_{L^{2}(\Omega_{s}(\omega))}\leq c\tau^{|\alpha|-3/2}e^{-\tau(s-t)a}\label{eq:2.5}
\end{equation}
for $|\alpha|\leq1$ and $s\geq t$, where $a>0$ is some constant
depending on $A_{t}(x')$. 

Without loss of generality, we consider the special case where $t=0$,
$\omega=e_{3}=(0,0,1)$ and choose $\eta=(1,0,0)$, $\zeta=(0,1,0)$.
The general case can be obtained from this special case by change
of coordinates. Define $L=L_{A}$ and $\widetilde{M}\cdot=e^{-i\tau x'\cdot\xi'}L(e^{i\tau x'\cdot\xi'}\cdot)$,
where $x'=(x_{1},x_{2})$ and $\xi'=(\xi_{1},\xi_{2})$ with $|\xi'|=1$,
then $\widetilde{M}$ is a differential operator. To be precise, by
using $a_{jl}=a_{lj}$, we calculate $\widetilde{M}$ to be given
by 
\begin{eqnarray*}
\widetilde{M} & = & -\tau^{2}\sum_{jl}a_{jl}\xi_{j}\xi_{l}+2\tau\sum_{jl}a_{jl}(i\xi_{l})\partial_{j}+\sum_{jl}a_{jl}\partial_{j}\partial_{l}\\
 &  & +\sum_{jl}(\partial_{j}a_{jl})(i\tau\xi_{l})+\sum_{jl}(\partial a_{jl})\partial_{l}+k^{2}\\
 & = & -\tau^{2}\sum_{jl}a_{jl}\xi_{j}\xi_{l}+2\tau\sum_{l}a_{3l}(i\xi_{l})\partial_{3}+a_{33}\partial_{3}\partial_{3}\\
 &  & +2\tau\sum_{j\neq3,l}a_{jl}(i\xi_{l})\partial_{j}+\sum_{(j,l)\backslash\{3,3\}}a_{jl}\partial_{j}\partial_{l}\\
 &  & +\sum_{jl}(\partial_{j}a_{jl})(i\tau\xi_{l})+\sum_{jl}(\partial_{j}a_{jl})\partial_{l}+k^{2}
\end{eqnarray*}
with $\xi_{3}=0$. Now, we want to solve 
\[
\widetilde{M}v=0,
\]
which is equivalent to $Mv=0$, where $M=a_{33}^{-1}\widetilde{M}$.
Now, we use the same idea in \cite{JN-ods}, define $\left\langle e,f\right\rangle =\sum_{ij}a_{ij}e_{i}f_{j}$,
where $e=(e_{1},e_{2},e_{3})$, $f=(f_{1},f_{2},f_{3})$ and denote
$\left\langle e,f\right\rangle _{0}=\left\langle e,f\right\rangle |_{x_{3}=0}$.
Let $P$ be a differential operator, and we define the order of $P$,
denoted by $ord(P)$, in the following sense:
\[
\|P(e^{-\tau x_{3}A(x')}\varphi(x')\|_{L^{2}(\mathbb{R}_{+}^{3})}\leq c\tau^{ord(P)-1/2},
\]
where $\mathbb{R}_{+}^{3}=\{x_{3}>0\}$, $A(x')$ is a smooth complex
function with its real part greater than 0 and $\varphi(x')\in C_{0}^{\infty}(\mathbb{R}^{2})$.
In this sense, similar to \cite{JN-ods}, we can see that $\tau$,
$\partial_{3}$ are of order 1, $\partial_{1},\partial_{2}$ are of
order 0 and $x_{3}$ is of order -1.

Now according to this order, the principal part $M_{2}$ (order 2)
of $M$ is:
\[
M_{2}=-\{D_{3}^{2}+2\tau\left\langle e_{3},e_{3}\right\rangle _{0}^{-1}\left\langle e_{3},\rho\right\rangle _{0}D_{3}+\tau^{2}\left\langle e_{3},e_{3}\right\rangle _{0}^{-1}\left\langle \rho,\rho\right\rangle _{0}\}
\]
with $D_{3}=-i\partial_{3}$ and $\rho=(\xi_{1},\xi_{2},0)$. Note
that the principal part $M_{2}$ does not involve the lower order
term $k^{2}\cdot$. Note that $M_{2}$ is obtained by the Taylor's
expansion of $M$ at $x_{3}=0$, that is, 
\begin{eqnarray*}
M(x',x_{3}) & = & M(x',0)+x_{3}\partial_{3}M(x',0)+\cdots+\dfrac{x_{3}^{N-1}}{(N-1)!}\partial_{3}^{N-1}M(x',0)+R\\
 & = & M_{2}+M_{1}+\cdots+M_{-N+1}+R,
\end{eqnarray*}
where ord$(M_{j})=j$ and ord$(R)=-N$. To solve $Mv=0$ is equivalent
to solve 
\begin{equation}
M_{2}v=-(M_{1}+\cdots+M_{-N+1}+R)v:=f.\label{eq:2.6}
\end{equation}
If we set $w_{1}=v$ and $w_{2}=-\tau^{-1}\left\langle e_{3},e_{3}\right\rangle _{0}D_{3}v-\left\langle e_{3},\rho\right\rangle _{0}v$,
then we can compute 
\begin{equation}
D_{3}w_{1}=-\tau\left\langle e_{3},e_{3}\right\rangle _{0}^{-1}\left\langle e_{3},\rho\right\rangle _{0}w_{1}-\tau\left\langle e_{3},e_{3}\right\rangle _{0}^{-1}w_{2}\label{eq:2.7}
\end{equation}
and 
\begin{eqnarray}
D_{3}w_{2} & = & -\tau\{\left\langle \rho,e_{3}\right\rangle _{0}^{2}\left\langle e_{3},e_{3}\right\rangle _{0}^{-1}-\left\langle \rho,\rho\right\rangle _{0}\}w_{1}-\tau\left\langle \rho,e_{3}\right\rangle _{0}\left\langle e_{3},e_{3}\right\rangle _{0}^{-1}w_{2}\label{eq:2.8}\\
 &  & +\tau^{-1}\left\langle e_{3},e_{3}\right\rangle _{0}f.\nonumber 
\end{eqnarray}
For detail calculations, we refer readers to see \cite{JN-ods}. If
we set $W=[w_{1},w_{2}]^{T}$ and use (\ref{eq:2.7}) and (\ref{eq:2.8}),
we have 
\[
D_{3}W=\tau KW+\left[\begin{array}{c}
0\\
\tau^{-1}\left\langle e_{3},e_{3}\right\rangle _{0}f
\end{array}\right],
\]
where 
\begin{equation}
K=\left[\begin{array}{cc}
\left\langle e_{3},e_{3}\right\rangle _{0}^{-1}\left\langle e_{3},\rho\right\rangle _{0} & \left\langle e_{3},e_{3}\right\rangle _{0}^{-1}\\
\left\langle \rho,e_{3}\right\rangle _{0}^{2}\left\langle e_{3},e_{3}\right\rangle _{0}^{-1}-\left\langle \rho,\rho\right\rangle _{0} & \left\langle \rho,e_{3}\right\rangle _{0}\left\langle e_{3},e_{3}\right\rangle _{0}^{-1}
\end{array}\right].\label{eq:2.9}
\end{equation}
By (\ref{eq:2.6}), we can express (\ref{eq:2.9}) as 
\begin{equation}
D_{3}W=(\tau K+K_{0}+\cdots+K_{-N}+S)W,\label{eq:2.10}
\end{equation}
where ord$(K_{j})=j$ and ord$(S)=-N-1$ and all the differential
operators $K_{j}$ involves only $x'$ derivatives. Moreover, $K$
is a matrix function independent of $x_{3}$ and its eigenvalues are
determined from 
\[
\det(\lambda I-K)=0,
\]
which is equivalent to 
\begin{equation}
\lambda^{2}-2\left\langle e_{3},e_{3}\right\rangle _{0}^{-1}\left\langle e_{3},\rho\right\rangle _{0}\lambda+\left\langle \rho,\rho\right\rangle _{0}=0.\label{eq:2.11}
\end{equation}
By using the uniform elliptic assumption on $(a_{ij})$ that (\ref{eq:2.11})
has roots $\lambda^{\pm}$ with Im$\lambda^{\pm}>0$. Similar to \cite{JN-ods},
we can set $\widetilde{Q}=[q^{+},q^{-}]$ be a nonsingular matrix
with linearly independent vectors $q^{\pm}$ such that 
\[
\widetilde{K}=\widetilde{Q}^{-1}K\widetilde{Q}=\left[\begin{array}{cc}
\lambda^{+} & 0\\
0 & \lambda^{-}
\end{array}\right],
\]
where $\lambda^{\pm}\in\mathbb{C}_{\pm}:=\{\pm\mbox{Im}\lambda>0\}$,
respectively. Moreover, we choose 
\begin{equation}
\widetilde{Q}=\left[\begin{array}{cc}
q & \overline{q}\\
q' & \overline{q'}
\end{array}\right],\label{eq:2.12}
\end{equation}
where 
\[
\left[\begin{array}{c}
q\\
q'
\end{array}\right]=[q^{+}]\mbox{ and }\left[\begin{array}{c}
\overline{q}\\
\overline{q'}
\end{array}\right]=[q^{-}]
\]
By virtue of the matrix $\widetilde{Q}$ in (\ref{eq:2.12}), we have
$\lambda^{-}=\overline{\lambda^{+}}$, and $\widetilde{Q}$ is nonsingular.
If we set $\widehat{Q}=\widetilde{Q}^{-1}W$, we get from (\ref{eq:2.10})
that 
\begin{equation}
D_{3}\widehat{W}=(\tau\widetilde{K}+\widehat{K}_{0}+\cdots+\widehat{K}_{-N}+S)\widehat{W},\label{eq:2.13}
\end{equation}
where ord$(\widehat{K}_{j})=j$ and ord($\widehat{S})=-N-1$. Similar
as before, we know that $\widehat{K}_{j}$ contains only $x'$ derivatives
since the original $K_{j}$ involves only $x'$ derivatives. In addition,
$\widehat{K}_{0}$ can be divided into terms involving $\tau x_{3}$
and terms formed by the differential operator in $\partial_{x'}$
with coefficients independent of $x_{3}$. Likewise, $\widehat{K}_{j}$
can be grouped into terms containing $\tau x_{3}^{-j+1},\tau^{-1}x_{3}^{j-1},x_{3}^{-j}$,
respectively, where $-N\leq j\leq-1$.

From now on, we have decoupled $K$ by choosing a suitable matrix
function $\widetilde{Q}$, next we want to decouple $\widehat{K}_{0},\cdots\widehat{K}_{-N}$.
First, we show how to decouple $\widehat{K}_{0}$. Let $\widehat{W}=(1+x_{3}A^{(0)}+\tau^{-1}B^{(0)})\widetilde{W}^{(0)}$
with $A^{(0)},B^{(0)}$ being differential operators in $\partial_{x'}$
with coefficients independent of $x_{3}$, then we have 
\begin{eqnarray*}
D_{3}\widehat{W}^{(0)} & = & \{\tau\widetilde{K}+(\widehat{K}_{0}-\tau x_{3}A^{(0)}\widetilde{K}+\tau x_{3}\widetilde{K}A^{(0)}-B^{(0)}\widetilde{K}+\widetilde{K}B^{(0)}+iA^{(0)})\\
 &  & +\widehat{K}'_{-1}+\cdots\}\widehat{W}^{(0)},
\end{eqnarray*}
where ord($\widehat{K}'_{-1})=-1$ and the remainder contains terms
of order at most -2. Let $\widetilde{K}_{0}:=\widehat{K}_{0}-\tau x_{3}A^{(0)}\widetilde{K}+\tau x_{3}\widetilde{K}A^{(0)}-B^{(0)}\widetilde{K}+\widetilde{K}B^{(0)}+iA^{(0)}$,
we analyze $\widetilde{K}_{0}$ more carefully. Set $\widehat{K}_{0}=\tau x_{3}\widehat{K}_{0,1}+\widehat{K}_{0,2}$
and express $\widehat{K}_{0,1},\widehat{K}_{0,2},A^{(0)}$ and $B^{(0)}$
in block forms, that is, 
\[
\widehat{K}_{0,l}=\left[\begin{array}{cc}
\widehat{K}_{0,l}(1,1) & \widehat{K}_{0,l}(1,2)\\
\widehat{K}_{0,l}(2,1) & \widehat{K}_{0,l}(2,2)
\end{array}\right],\mbox{ }l=1,2,
\]
\[
A^{(0)}=\left[\begin{array}{cc}
A^{(0)}(1,1) & A^{(0)}(1,2)\\
A^{(0)}(2,1) & A^{(0)}(2,2)
\end{array}\right]\mbox{ and }B^{(0)}=\left[\begin{array}{cc}
B^{(0)}(1,1) & B^{(0)}(1,2)\\
B^{(0)}(2,1) & B^{(0)}(2,2)
\end{array}\right].
\]
Then the off-diagonal blocks of $\widetilde{K}_{0}$ are given by:
\begin{eqnarray*}
\widetilde{K}_{0}(1,2) & = & \tau x_{3}\{\widehat{K}_{0,1}(1,2)-A^{(0)}(1,2)\lambda^{-}+\lambda^{+}A^{(0)}(1,2)\}\\
 &  & +\{\widehat{K}_{0,2}(1,2)+iA^{(0)}(1,2)-B^{(0)}(1,2)\lambda^{-}+\lambda^{+}B^{(0)}(1,2)\},
\end{eqnarray*}
\begin{eqnarray*}
\widetilde{K}_{0}(2,1) & = & \tau x_{3}\{\widehat{K}_{0,1}(2,1)-A^{(0)}(2,1)\lambda^{-}+\lambda^{+}A^{(0)}(2,1)\}\\
 &  & +\{\widehat{K}_{0,2}(2,1)+iA^{(0)}(2,1)-B^{(0)}(2,1)\lambda^{-}+\lambda^{+}B^{(0)}(2,1)\}.
\end{eqnarray*}
Since $\lambda^{\pm}\in\mathbb{C}_{\pm}$, we can find suitable $A^{(0)}(1,2)$
and $A^{(0)}(2,1)$ such that 
\[
\begin{cases}
\widehat{K}_{0,1}(1,2)-A^{(0)}(1,2)\lambda^{-}+\lambda^{+}A^{(0)}(1,2)=0\\
\widehat{K}_{0,1}(2,1)-A^{(0)}(2,1)\lambda^{-}+\lambda^{+}A^{(0)}(2,1)=0
\end{cases}
\]
(see similar arguments in \cite{Taylor}). Similarly, we can use the
same method to find $B^{(0)}(1,2)$ and $B^{(0)}(2,1)$ so that 
\begin{equation}
\begin{cases}
\widehat{K}_{0,2}(1,2)+iA^{(0)}(1,2)-B^{(0)}(1,2)\lambda^{-}+\lambda^{+}B^{(0)}(1,2)=0,\\
\widehat{K}_{0,2}(2,1)+iA^{(0)}(2,1)-B^{(0)}(2,1)\lambda^{-}+\lambda^{+}B^{(0)}(2,1)=0.
\end{cases}\label{eq:2.14}
\end{equation}
Since $\widehat{K}_{0,2}(1,2)$ and $\widehat{K}_{0,2}(2,1)$ are
differential operators in $\partial_{x'}$ with coefficients independent
of $x_{3}$, we will look for $B^{(0)}(1,2)$ and $B^{(0)}(2,1)$
as the same type of differential operators. By (\ref{eq:2.14}) and
using $\lambda^{\pm}\in\mathbb{C}_{\pm}$, we can solve for $B^{(0)}(1,2)$
and $B^{(0)}(2,1)$. To find $A^{(0)}$ and $B^{(0)}$, we simply
set diagonal blocks of them are zero, i.e., 
\[
A^{(0)}=\left[\begin{array}{cc}
0 & A^{(0)}(1,2)\\
A^{(0)}(2,1) & 0
\end{array}\right]\mbox{ and }B^{(0)}=\left[\begin{array}{cc}
0 & B^{(0)}(1,2)\\
B^{(0)}(2,1) & 0
\end{array}\right].
\]
With these matrices $A^{(0)}$ and $B^{(0)}$, we can see that 
\begin{equation}
D_{3}\widetilde{W}^{(0)}=\{\tau\widetilde{K}+\widetilde{K}_{0}+\widehat{K}'_{-1}+\cdots\}\widetilde{W}^{(0)}\label{eq:2.15}
\end{equation}
where 
\[
\widetilde{K}_{0}=\left[\begin{array}{cc}
\widetilde{K}_{0}(1,1) & 0\\
0 & \widetilde{K}_{0}(2,2)
\end{array}\right].
\]

Moreover, we want to decouple $\widehat{K}'_{-1}$ and $\widehat{K}'_{-1}$
can be written as $\widehat{K}'_{-1}=\tau x_{3}^{2}\widehat{K}'_{-1,1}+x_{3}\widehat{K}'_{-1,2}+\tau^{-1}\widehat{K}'_{-1,3}$.
We can see that $\widehat{K}'_{-1.1}$, $\widehat{K}'_{-1,2}$ and
$\widehat{K}'_{-1,3}$ are differential operators in $\partial_{x'}$
of order zero, one and two with coefficients independent of $x_{3}$,
respectively. Similarly, we can set $\widetilde{W}^{(0)}=(I+x_{3}^{2}A^{(1)}+\tau^{-1}x_{3}B^{(1)}+\tau^{-2}C^{(1)})\widetilde{W}^{(1)}$,
where $A^{(1)}$, $B^{(1)}$ and $C^{(1)}$ are differential operators
in $\partial_{x'}$. Now plugging $\widetilde{W}^{(0)}$ of above
form into (\ref{eq:2.15}), we have 
\begin{eqnarray}
D_{3}\widetilde{W}^{(1)} & = & \{\tau\widetilde{K}+\widetilde{K}_{0}+\tau x_{3}^{2}(\widehat{K}'_{-1,1}-A^{(1)}\widetilde{K}+\widetilde{K}A^{(1)})+x_{3}(\widehat{K}'_{-1,2}-B^{(1)}\widetilde{K}\nonumber \\
 &  & +\widetilde{K}B^{(1)}+2A^{(1)})+\tau^{-1}(\widehat{K}'_{-1,3}-C^{(1)}\widetilde{K}+\widetilde{K}C^{(1)}+iB^{(1)})]\nonumber \\
 &  & +\cdots\}\widetilde{W}^{(1)}\label{eq:2.16}
\end{eqnarray}
where the remainder consists of terms with order at most -2. Then
we use the same argument, we can find suitable $A^{(1)}$, $B^{(1)}$
and $C^{(1)}$ such that the off-diagonal blocks of the order -1 term
on the right hand side of (\ref{eq:2.16}) are zero. Therefore, we
obtain 
\[
D_{3}\widetilde{W}^{(1)}=\{\tau\widetilde{K}+\widetilde{K}_{0}+\widetilde{K}_{-1}+\cdots\}\widetilde{W}^{(1)}
\]
with 
\[
\widetilde{K}_{-1}=\left[\begin{array}{cc}
\widetilde{K}_{-1}(1,1) & 0\\
0 & \widetilde{K}_{-1}(2,2)
\end{array}\right].
\]
Recursively, by defining 
\begin{eqnarray*}
\widehat{W} & = & (I+x_{3}A^{(0)}+\tau^{-1}B^{(0)})(I+x_{3}^{2}A^{(1)}+\tau^{-1}x_{3}B^{(1)}+\tau^{-2}C^{(1)})\cdots\\
 &  & (I+x_{3}^{N+1}A^{(N)}+\tau^{-1}x_{3}^{N}B^{(N)}+\tau^{-2}x_{3}^{N-1}C^{(N)})\widetilde{W}^{(N)}
\end{eqnarray*}
with suitable $A^{(j)}$, $B^{(j)}$ and $C^{(j)}$ for $0\leq j\leq N$
($C^{(0)}=0$), we can transform the equation (\ref{eq:2.13}) into
\begin{equation}
D_{3}\widetilde{W}^{(N)}=\{\tau\widetilde{K}+\widetilde{K}_{0}+\cdots+\widetilde{K}_{-N}+\widetilde{S}\}\widetilde{W}^{(N)},\label{eq:2.17}
\end{equation}
where $\widetilde{K}_{-j}$ for all $0\leq j\leq N$ are decoupled
and ord$(\widetilde{S})=-N-1$. Note that all diagonal blocks of $A^{(j)}$
and $B^{(j)}$ are zero.

Now in view of (\ref{eq:2.17}), we consider the equation
\[
D_{3}\hat{v}^{(N)}=\{\tau\lambda^{+}+\widetilde{K}_{0}(1,1)+\cdots+\widetilde{K}_{-N}(1,1)\}\hat{v}^{(N)},
\]
with an approximated solution of the form 
\[
\hat{v}^{(N)}=\sum_{j=0}^{N+1}\hat{v}_{-j}^{(N)},
\]
where $\hat{v}_{-j}^{(N)}$ for $0\leq j\leq N$ satisfy 
\[
\begin{cases}
D_{3}\hat{v}_{0}^{(N)}=\tau\lambda^{+}\hat{v}_{0}^{(N)}, & \hat{v}_{0}^{(N)}|_{x_{3}=0}=\chi_{t}(x')b\\
D_{3}\hat{v}_{-1}^{(N)}=\tau\lambda^{+}\hat{v}_{-1}^{(N)}+\widetilde{K}_{0}(1,1)\hat{v}_{0}^{(N)}, & \hat{v}_{-1}^{(N)}|_{x_{3}=0}=0\\
\vdots & \vdots\\
D_{3}\hat{v}_{-N-1}^{(N)}=\tau\lambda^{+}\hat{v}_{-N-1}^{(N)}+\sum_{j=0}^{N}\widetilde{K}_{-j}(1,1)\hat{v}_{j-N}^{(N)}, & \hat{v}_{-N-1}^{(N)}|_{x_{3}=0}=0,
\end{cases}
\]
where $\chi_{t}(x')\in C_{0}^{\infty}(\mathbb{R}^{2})$ and $b\in\mathbb{C}$.
It easy to solve $\hat{v}_{0}^{(N)}=\exp(i\tau x_{3}\lambda^{+})\chi_{t}(x')b$
and $\hat{v}_{-1}^{(N)}=\exp(i\tau x_{3}\lambda^{+})\int_{0}^{x_{3}}\exp(-i\tau s\lambda^{+})\widetilde{K}_{0}(1,1)\hat{v}_{0}^{(N)}ds$.
Moreover, we can use the ord($x_{3})=-1$ and ord$(\partial_{j})=0$
with $j=1,2$ to derive that 
\[
\|x_{3}^{\beta}\partial_{x'}^{\alpha}\hat{v}_{0}^{(N)}\|_{L^{2}(\mathbb{R}_{+}^{3})}\leq c\tau^{-\beta-1/2}
\]
for $\beta\in\mathbb{Z}_{+}$ and multi-index $\alpha$. Similarly,
we can compute
\begin{equation}
\|\hat{v}_{-1}^{(N)}\|_{L^{2}(\mathbb{R}_{+}^{3})}^{2}\leq c\tau^{-3}.\label{eq:2.18}
\end{equation}
For the derivation of (\ref{eq:2.18}), it can be found in \cite{JN-ods}.
Moreover, by similar computations we can show that 
\[
\|x_{3}^{\beta}\partial_{x'}^{\alpha}(\hat{v}_{-1}^{(N)})\|_{L^{2}(\mathbb{R}_{+}^{3})}\leq c\tau^{-\beta-3/2}
\]
and for $\hat{v}_{-j}^{(N)}$, $j=2,\ldots,N+1$, we have 
\[
\|x_{3}^{\beta}\partial_{x'}^{\alpha}(\hat{v}_{-j}^{(N)})\|_{L^{2}(\mathbb{R}_{+}^{3})}\leq c\tau^{-\beta-j-1/2}
\]
for $2\leq j\leq N+1$.

Thus, if we set $V^{(N)}=\left[\begin{array}{c}
\hat{v}^{(N)}\\
0
\end{array}\right]$, then we have 
\[
\begin{cases}
D_{3}V^{(N)}-\{\tau\widetilde{K}+\widetilde{K}_{0}+\cdots+\widetilde{K}_{-N}\}V^{(N)}=\widetilde{R},\\
V^{(N)}|_{x_{3}=0}=\left[\begin{array}{c}
\chi_{t}(x')b\\
0
\end{array}\right],
\end{cases}
\]
where 
\[
\|\widetilde{R}\|_{L^{2}(\mathbb{R}_{+}^{3})}\leq c\tau^{-N-3/2}.
\]
Define $v$ to be the function of the first component of $\widetilde{Q}(I+x_{3}A^{(0)}+\tau^{-1}B^{(0)})(I+x_{3}^{2}A^{(1)}+\tau^{-1}x_{3}B^{(1)}+\tau^{-2}C^{(1)})\cdots(I+x_{3}^{N+1}A^{(N)}+\tau^{-1}x_{3}^{N}B^{(N)}+\tau^{-2}x_{3}^{N-1}C^{(N)})V^{(N)}$
and set $w=\exp(i\tau x'\cdot\xi')\tilde{v}$, we have 
\begin{eqnarray*}
w & = & q\exp(i\tau x'\cdot\xi')\exp(i\tau x_{3}\lambda^{+}(x'))\chi_{t}(x')b+\exp(i\tau x'\cdot\xi')\tilde{\gamma}(x,\tau)\\
 &  & q\exp(i\tau x'\cdot\xi')\exp(-\tau x_{3}(-i\lambda^{+}(x')))\chi_{t}(x')b+\gamma(x,\tau)
\end{eqnarray*}
and 
\[
w|_{x_{3}=0}=\exp(i\tau x'\cdot\xi')\{\chi_{t}(x')qb+\beta_{0}(x',\tau)\},
\]
where $\gamma$ satisfies the estimate (\ref{eq:2.5}) on $\Omega_{s}:=\{x_{3}>s\}\cap\Omega$
for $s\ge0$ and $\beta_{0}(x',\tau)=\tilde{\gamma}(x',0,\tau)$ is
supported in supp($\chi_{t})$ with $\|\beta_{0}(\cdot,\tau)\|_{L^{\infty}}\leq c\tau^{-1}$.
Also, we have 
\[
\|M\tilde{v}\|_{L^{2}(\Omega_{0})}\leq c\tau^{-N-1/2}.
\]

Let $u=w+r=e^{ix'\cdot\xi'}\tilde{v}+r$ and r be the solution to
the boundary value problem
\begin{equation}
\begin{cases}
Lr=-e^{ix'\cdot\xi'}\widetilde{M}\tilde{v} & \mbox{ in }\Omega_{0},\\
r=0 & \mbox{ on }\partial\Omega_{0}.
\end{cases}\label{eq:2.19}
\end{equation}
The existence of $r$ solving (\ref{eq:2.19}) is by using the Lax-Milgram
theorem and we have the following estimate 
\[
\|r\|_{H^{1}(\Omega_{0})}\leq c\tau^{-N-1/2},
\]
which is the estimate (\ref{eq:2.4}) on $\Omega_{0}$. We complete
the construction of the oscillating-decaying solutions for the case
$t=0$ and $\omega=(0,0,1)$ in the anisotropic elliptic equations
case. The oscillating-decaying solution in the general case can be
obtained by using change of coordinates.

\section{Tools and estimates}

In this section, we introduce the Runge approximation property and
a very useful elliptic estimate: Meyers $L^{p}$-estimates.

\subsection{Runge approximation property}
\begin{defn}
\cite{Lax} Let $L$ be a second order elliptic operator, solutions
of an equation $Lu=0$ are said to have the Runge approximation property
if, whenever $K$ and $\Omega$ are two simply connected domains with
$K\subset\Omega$, any solution in $K$ can be approximated uniformly
in compact subsets of $K$ by a sequence of solutions which can be
extended as solution to $\Omega$.
\end{defn}
There are many applications for Runge approximation property in inverse
problems. Similar results for some elliptic operators can be found
in \cite{Lax}, \cite{Malgrange}. The following theorem is a classical
result for Runge approximation property for a second order elliptic
equation.
\begin{thm}
(Runge approximation property) Let $L_{0}\cdot=\nabla(A^{0}(x)\nabla\cdot)+k^{2}\cdot$
be a second order elliptic differential operator with $A^{0}(x)$
to be Lipschitz. Assume that $k^{2}$ is not a Dirichlet eigenvalue
of $-\nabla(A^{0}(x)\nabla\cdot)$. Let $O$ and $\Omega$ be two
open bounded domains with smooth boundary in $\mathbb{R}^{3}$ such
that $O$ is convex and $\bar{O}\subset\Omega$. 

Let $u_{0}\in H^{1}(O)$ satisfy 
\[
L_{0}u_{0}=0\mbox{ in }O.
\]
Then for any compact subset $K\subset O$ and any $\epsilon>0$, there
exists $U\in H^{1}(\Omega)$ satisfying 
\[
L_{0}U=0\mbox{ in }\Omega,
\]
such that 
\[
\|u_{0}-U\|_{H^{1}(K)}\leq\epsilon.
\]

\end{thm}
Note that we have assumed that $A^{0}\in B^{\infty}(\mathbb{R}^{3})$,
it is easy to see $A^{0}(x)$ is a Lipschitz continuous function,
it possesses the weak continuation property.The proof can be found
in \cite{Lax} and \cite{JN-ods}, we omit details here.

\subsection{Elliptic estimates and some identities}

We need some estimates for solutions to some Dirichlet problems which
will be used in next section. Recall that, for $f\in H^{1/2}(\partial\Omega)$,
let $u$ and $u_{0}$ be solutions to the Dirichlet problems (\ref{eq:1.1})
and (\ref{eq:1.2}), respectively. Note that $a_{ij}(x)=a_{ij}^{0}(x)\chi_{\Omega\backslash D}+\widetilde{a_{ij}}(x)\chi_{D}$
and we set $w=u-u_{0}$, then $w$ satisfies the Dirichlet problem
\begin{equation}
\begin{cases}
\nabla\cdot(A(x)\nabla w)+k^{2}w=-\nabla\cdot((\widetilde{A}\chi_{D}-A^{0}\chi_{D})\nabla u_{0}) & \mbox{ in }\Omega\\
w=0 & \mbox{ on }\partial\Omega
\end{cases}\label{eq:3.1}
\end{equation}
where $A(x)=(a_{ij}(x))$, $A^{0}(x)=(a_{ij}^{0}(x))$ and $\widetilde{A}(x)=(\widetilde{a_{ij}}(x))$.
Then we have some estimates for $w$.
\begin{lem}
There exists a positive constant $C$ independent of $w$ such that
we have 
\[
\|w\|_{L^{2}(\Omega)}\leq C\|\nabla w\|_{L^{p}(\Omega)}
\]
for $\dfrac{6}{5}\leq p\leq2$ if $n=3$.
\end{lem}
The proof follow from \cite{Sini} by Freidrichs inequality, see \cite{Majza}
p.258 and use a standard elliptic regularity.
\begin{lem}
There exists $\epsilon\in(0,1)$, depending only on $\Omega$, $A^{0}(x)=(a_{ij}^{0}(x))$
and $\widetilde{A}(x)=(\widetilde{a_{ij}}(x))$ such that 
\[
\|\nabla w\|_{L^{p}(\Omega)}\leq C\|u_{0}\|_{W^{1,p}(D)}
\]
for $\max\{2-\epsilon,\dfrac{6}{5}\}<p\leq2$ if $n=3$.\end{lem}
\begin{proof}
The proof is also followed from \cite{Sini}. Set $f:=-(\widetilde{A}\chi_{D}-A^{0}\chi_{D})\nabla u_{0}$,
$h:=0$. Let $w_{0}$ be a solution of 
\begin{equation}
\begin{cases}
\nabla\cdot(A(x)\nabla w_{0})+k^{2}w_{0}=\nabla\cdot f & \mbox{ in }\Omega,\\
w_{0}=0 & \mbox{ on }\partial\Omega.
\end{cases}\label{eq:3.2}
\end{equation}
The following $L^{p}$-estimate of $w_{0}$, followed from \cite{Meyer},
then we can get 
\begin{equation}
\|\nabla w_{0}\|_{L^{p}(\Omega)}\leq C\|f\|_{L^{p}(\Omega)}\label{eq:3.3}
\end{equation}
for $p\in(\max\{2-\epsilon,\dfrac{6}{5}\},2]$, where $\epsilon\in(0,1)$
depends on $\Omega$, $A^{0}(x)=(a_{ij}^{0}(x))$ and $\widetilde{A}(x)=(\widetilde{a_{ij}}(x))$.
We set $W:=w-w_{0}$, then since $w=w_{0}+W$, we have 
\begin{equation}
\|\nabla w\|_{L^{p}(\Omega)}\leq C(\|\nabla w_{0}\|_{L^{p}(\Omega)}+\|\nabla W\|_{L^{p}(\Omega)}).\label{eq:3.4}
\end{equation}
Moreover, $W$ satisfies 
\begin{equation}
\begin{cases}
\nabla\cdot(A(x)\nabla W)+k^{2}W=0 & \mbox{ in }\Omega,\\
W=0 & \mbox{ on }\partial\Omega.
\end{cases}\label{eq:3.5}
\end{equation}
By the standard elliptic regularity, we have 
\[
\|W\|_{H^{1}(\Omega)}\leq C\|w_{0}\|_{L^{2}(\Omega)}.
\]
Thus, we get for $p\leq2$, 
\begin{equation}
\|\nabla W\|_{L^{p}(\Omega)}\leq C\|\nabla W\|_{L^{2}(\Omega)}\leq C\|W\|_{H^{1}(\Omega)}\leq C\|w_{0}\|_{L^{2}(\Omega)}.\label{eq:3.6}
\end{equation}

By Sobolev embedding theorem, we get 
\begin{equation}
\|w_{0}\|_{L^{2}(\Omega)}\leq C\|w_{0}\|_{W^{1.p}(\Omega)}\label{eq:3.7}
\end{equation}
for $p\geq\dfrac{6}{5}$ if $n=3$. Use Poincar$\acute{e}$'s inequality
in $L^{p}$ spaces ($w_{0}|_{\partial\Omega}=0$), we have 
\begin{equation}
\|w_{0}\|_{L^{2}(\Omega)}\leq C\|\nabla w_{0}\|_{L^{p}(\Omega)}\label{eq:3.8}
\end{equation}
for $p\geq\dfrac{6}{5}$ if $n=3$. Combining (\ref{eq:3.3}) with
(\ref{eq:3.4}), (\ref{eq:3.6}) and (\ref{eq:3.8}), we can obtain
\[
\|\nabla w\|_{L^{p}(\Omega)}\leq C\|f\|_{L^{p}(\Omega)}\leq C\|u_{0}\|_{W^{1,p}(D)}
\]
for $\max\{2-\epsilon,\dfrac{6}{5}\}<p\leq2$ if $n=3$.
\end{proof}
Recall the Dirichlet-to-Neumann map which we have defined in the section
1: $\Lambda_{D}f:=A\nabla u\cdot\nu$ and $\Lambda_{\emptyset}f:=A^{0}\nabla u_{0}\cdot\nu$,
where $\nu=(\nu_{1},\nu_{2},\nu_{3})$ is an outer normal on $\partial\Omega$.
We next prove some useful identities.
\begin{lem}
$\int_{\partial\Omega}(\Lambda_{D}-\Lambda_{\emptyset})f\bar{f}d\sigma=\mbox{Re}\int_{D}(\widetilde{A}-A^{0})\nabla u_{0}\cdot\overline{\nabla u}dx$.\end{lem}
\begin{proof}
It is clear that 
\begin{eqnarray*}
\int_{\partial\Omega}A\nabla u\cdot\nu\bar{\varphi}d\sigma & = & \int_{\Omega}\nabla\cdot(A\nabla u\bar{\varphi})dx\\
 & = & \int_{\Omega}\nabla\cdot(A\nabla u)\bar{\varphi}+A\nabla u\cdot\overline{\nabla\varphi}dx\\
 & = & -k^{2}\int_{\Omega}u\bar{\varphi}dx+\int_{\Omega}A\nabla u\cdot\overline{\nabla\varphi}dx
\end{eqnarray*}
$\forall\varphi\in H^{1}(\Omega)$. Since $u=u_{0}=f$ on $\partial\Omega$,
the left hand side of the identity has the same value whether we take
$\varphi=u$ or $\varphi=u_{0}$, and it is equal to $\int_{\partial\Omega}\Lambda_{D}f\bar{f}d\sigma$.
\begin{eqnarray*}
\int_{\partial\Omega}\Lambda_{D}f\bar{f}d\sigma & = & -k^{2}\int_{\Omega}u\overline{u_{0}}dx+\int_{\Omega}A\nabla u\cdot\overline{\nabla u_{0}}dx\\
 & = & -k^{2}\int_{\Omega}|u|^{2}dx+\int_{\Omega}A\nabla u\cdot\overline{\nabla u}dx.
\end{eqnarray*}
The right hand side of the identity above is real. Hence, by taking
the real part, we have 
\[
\int_{\partial\Omega}\Lambda_{D}f\bar{f}d\sigma=-k^{2}\mbox{Re}\int_{\Omega}u\overline{u_{0}}dx+\mbox{Re}\int_{\Omega}A\nabla u\cdot\overline{\nabla u_{0}}dx
\]
and 
\[
\int_{\partial\Omega}\Lambda_{\emptyset}f\bar{f}d\sigma=-k^{2}\mbox{Re}\int_{\Omega}u\overline{u_{0}}dx+\mbox{Re}\int_{\Omega}A^{0}\nabla u\cdot\overline{\nabla u_{0}}dx.
\]
Therefore, we have 
\begin{eqnarray}
\int_{\partial\Omega}(\Lambda_{D}-\Lambda_{\emptyset})f\bar{f}d\sigma & = & \mbox{Re}\int_{\Omega}(A-A^{0})\nabla u\cdot\overline{\nabla u_{0}}dx\label{eq:3.9}\\
 & = & \mbox{Re}\int_{\Omega}(\widetilde{A}-A^{0})\chi_{D}\nabla u\cdot\overline{\nabla u_{0}}dx.\nonumber 
\end{eqnarray}

\end{proof}
The estimates in the following lemma play an important role in our
reconstruction algorithm.
\begin{lem}
We have the following identities:
\begin{eqnarray}
\int_{\partial\Omega}(\Lambda_{D}-\Lambda_{\emptyset})f\bar{f}d\sigma & = & -\int_{\Omega}A\nabla w\cdot\overline{\nabla w}dx+k^{2}\int_{\Omega}|w|^{2}dx\label{eq:3.10}\\
 &  & +\int_{D}(A^{0}-\widetilde{A})\nabla u_{0}\cdot\overline{\nabla u_{0}}dx,\nonumber 
\end{eqnarray}
\begin{eqnarray}
\int_{\partial\Omega}(\Lambda_{D}-\Lambda_{\emptyset})f\bar{f}d\sigma & = & \int_{\Omega}A^{0}\nabla w\cdot\overline{\nabla w}dx-k^{2}\int_{\Omega}|w|^{2}dx\label{eq:3.11}\\
 &  & +\int_{D}(\widetilde{A}-A^{0})\nabla u\cdot\overline{\nabla u}dx.\nonumber 
\end{eqnarray}

In particular, we have 
\begin{equation}
\int_{\partial\Omega}(\Lambda_{D}-\Lambda_{\emptyset})f\bar{f}d\sigma\leq k^{2}\int_{\Omega}|w|^{2}dx+\widehat{\Lambda}\int_{D}|\nabla u_{0}|^{2}dx,\label{eq:3.12}
\end{equation}
\begin{equation}
\int_{\partial\Omega}(\Lambda_{D}-\Lambda_{\emptyset})f\bar{f}d\sigma\geq c\int_{\Omega}|\nabla u_{0}|^{2}dx-k^{2}\int_{\Omega}|w|^{2}dx,\label{eq:3.13}
\end{equation}
where $c$ depending only on $\widetilde{\lambda}$ and $\lambda^{0}$.\end{lem}
\begin{proof}
Multiplying the identity 
\[
\nabla\cdot(A(x)\nabla w)+k^{2}w+\nabla\cdot((\widetilde{A}\chi_{D}-A^{0}\chi_{D})\nabla u_{0})=0
\]
by $\bar{w}$ and integrating over $\Omega$, we get 
\begin{eqnarray*}
0 & = & \int_{\Omega}\nabla\cdot(A\nabla w)\bar{w}dx+\int_{\Omega}\nabla\cdot((A^{0}-\widetilde{A})\chi_{D}\nabla u_{0})\bar{w}dx+k^{2}\int_{\Omega}|w|^{2}dx\\
 & = & -\int_{\Omega}A\nabla w\cdot\overline{\nabla w}dx+\int_{\partial\Omega}A\dfrac{\partial w}{\partial\nu}\bar{w}d\sigma-\int_{\Omega}(A^{0}-\widetilde{A})\chi_{D}\nabla u_{0}\cdot\overline{\nabla w}dx\\
 &  & +\int_{\partial\Omega}(A^{0}-\widetilde{A})\chi_{D}\dfrac{\partial u_{0}}{\partial\nu}\bar{w}d\sigma+k^{2}\int_{\Omega}|w|^{2}dx\\
 & = & -\int_{\Omega}A\nabla w\cdot\overline{\nabla w}dx-\int_{D}(A^{0}-\widetilde{A})\nabla u_{0}\cdot\overline{\nabla w}dx+k^{2}\int_{\Omega}|w|^{2}dx\\
 & = & \int_{\Omega}A\nabla w\cdot\overline{\nabla w}dx-\int_{D}(A^{0}-\widetilde{A})\nabla u_{0}\cdot\overline{\nabla u}dx+k^{2}\int_{\Omega}|w|^{2}dx\\
 &  & +\int_{\Omega}(A^{0}-\widetilde{A})\chi_{D}\nabla u_{0}\cdot\overline{\nabla u_{0}}dx,
\end{eqnarray*}
and use (\ref{eq:3.9}) we can obtain 
\[
\int_{\partial\Omega}(\Lambda_{D}-\Lambda_{\emptyset})f\bar{f}d\sigma=-\int_{\Omega}A\nabla w\cdot\overline{\nabla w}dx+\int_{D}(A^{0}-\widetilde{A}){}_{D}\nabla u_{0}\cdot\overline{\nabla u_{0}}dx+k^{2}\int_{\Omega}|w|^{2}dx.
\]

Similarly, multiplying the identity 
\[
0=\nabla\cdot((\widetilde{A}-A^{0})\chi_{D}\nabla u)+\nabla\cdot(A^{0}\nabla w)+k^{2}w=0
\]
by $\bar{w}$ and integrating over $\Omega$, we get 
\begin{eqnarray*}
0 & = & \int_{\Omega}\nabla\cdot((\widetilde{A}-A^{0})\chi_{D}\nabla u)\bar{w}dx+\int_{\Omega}\nabla\cdot(A^{0}\nabla w)\bar{w}dx+k^{2}\int_{\Omega}|w|^{2}dx\\
 & = & -\int_{D}(\widetilde{A}-A^{0})\nabla u\cdot\overline{\nabla w}dx-\int_{\Omega}A^{0}\nabla w\cdot\overline{\nabla w}dx+k^{2}\int_{\Omega}|w|^{2}dx\\
 & = & -\int_{D}(\widetilde{A}-A^{0})\nabla u\cdot\overline{\nabla u}dx-\int_{D}(\widetilde{A}-A^{0})\nabla u\cdot\overline{\nabla u_{0}}dx+k^{2}\int_{\Omega}|w|^{2}dx\\
 &  & -\int_{\Omega}A^{0}\nabla w\cdot\overline{\nabla w}dx,
\end{eqnarray*}
and use (\ref{eq:3.9}) again, we can obtain
\[
\int_{\partial\Omega}(\Lambda_{D}-\Lambda_{\emptyset})f\bar{f}d\sigma=\int_{\Omega}A^{0}\nabla w\cdot\overline{\nabla w}dx-k^{2}\int_{\Omega}|w|^{2}dx+\int_{D}(\widetilde{A}-A^{0})\nabla u\cdot\overline{\nabla u}dx.
\]

For the remaining part, (\ref{eq:3.12}) is an easy consequence of
(\ref{eq:3.10}) 
\begin{eqnarray*}
\int_{\partial\Omega}(\Lambda_{D}-\Lambda_{\emptyset})f\bar{f}d\sigma & \leq & k^{2}\int_{\Omega}|w|^{2}dx+\int_{D}(A^{0}-\widetilde{A})\nabla u_{0}\cdot\overline{\nabla u_{0}}dx\\
 & \leq & k^{2}\int_{\Omega}|w|^{2}dx+\widehat{\Lambda}\int_{D}|\nabla u_{0}|^{2}dx
\end{eqnarray*}

Finally, for the lower bound, we use 
\begin{eqnarray*}
A^{0}\nabla w\cdot\overline{\nabla w}+(\widetilde{A}-A^{0})\nabla u\cdot\overline{\nabla u} & = & \widetilde{A}\nabla u\cdot\overline{\nabla u}-2\mbox{Re}A^{0}\nabla u\cdot\overline{\nabla u_{0}}+A^{0}\nabla u_{0}\cdot\overline{\nabla u_{0}}\\
 & = & \widetilde{A}(\nabla u-(\widetilde{A})^{-1}A^{0}\nabla u_{0})\cdot(\overline{\nabla u-(\widetilde{A})^{-1}A^{0}\nabla u_{0}})\\
 &  & +(A^{0}-(\widetilde{A})^{-1}(A^{0})^{2})\nabla u_{0}\cdot\overline{\nabla u_{0}}\\
 & \geq & (A^{0}-(\widetilde{A})^{-1}(A^{0})^{2})\nabla u_{0}\cdot\overline{\nabla u_{0}}\\
 & \geq & c|\nabla u_{0}|^{2},
\end{eqnarray*}
since $\widetilde{A}(\nabla u-(\widetilde{A})^{-1}A^{0}\nabla u_{0})\cdot(\overline{\nabla u-(\widetilde{A})^{-1}A^{0}\nabla u_{0}})\geq0$
and note that $A^{0}-(\widetilde{A})^{-1}(A^{0})^{2}=(\widetilde{A})^{-1}(\widetilde{A}-A^{0})A^{0}$
has a positive lower bound depending only on $\widetilde{\lambda}$
and $\lambda^{0}$.
\end{proof}
Before stating our main theorem, we need to estimate $\|w\|_{L^{2}(\Omega)}$.
Fortunately, we can use Meyers $L^{p}$ estimates to help us to overcome
the difficulties (see lemma 3.2 and lemma 3.3). For the upper bound
of $\int_{\partial\Omega}(\Lambda_{D}-\Lambda_{\emptyset})f\bar{f}d\sigma$,
see (\ref{eq:3.11}), we use $\|w\|_{L^{2}(\Omega)}\leq C\|u_{0}\|_{W^{1,p}(D)}$
for $p\leq2$. Then we have 
\begin{equation}
\int_{\partial\Omega}(\Lambda_{D}-\Lambda_{\emptyset})f\bar{f}d\sigma\leq C\|u_{0}\|_{W^{1,p}(D)}^{2}.\label{eq:3.14}
\end{equation}

By (\ref{eq:3.13}) and the Meyers $L^{p}$ estimate $\|w\|_{L^{2}(\Omega)}\leq C\|u_{0}\|_{W^{1,p}(D)}$,
we have 
\begin{equation}
\int_{\partial\Omega}(\Lambda_{D}-\Lambda_{\emptyset})f\bar{f}d\sigma\ge c\int_{\Omega}|\nabla u_{0}|^{2}dx-c\|u_{0}\|_{W^{1,p}(D)}^{2}.\label{eq:3.15}
\end{equation}

\section{Detecting the convex hull of the unknown obstacle}

\subsection{Main theorem}

Recall that we have constructed the oscillating-decaying solutions
in section 2, and note that this solution can not be defined on the
whole domain, that is, the oscillating-decaying solutions $u_{\chi_{t},b,t,N,\omega}(x,\tau)$
only defined on $\Omega_{t}(\omega)\subsetneq\Omega$. Nevertheless,
with the help of the Runge approximation property, we can prove that
one can determine the convex hull of the unknown obstacle $D$ by$\Lambda_{D}f$
for infinitely many $f$.

We define $B$ to be an open ball in $\mathbb{R}^{3}$ such that $\overline{\Omega}\subset B$.
Assume that $\widetilde{\Omega}\subset\mathbb{R}^{3}$ is an open
Lipschitz domain with $\overline{B}\subset\widetilde{\Omega}$. As
in the section 2, set $\omega\in S^{2}$ and $\{\eta,\zeta,\omega\}$
forms an orthonormal basis of $\mathbb{R}^{3}$. Suppose $t_{0}=\inf_{x\in D}x\cdot\omega=x_{0}\cdot\omega$,
where $x_{0}=x_{0}(\omega)\in\partial D$. For any $t\leq t_{0}$
and $\epsilon>0$ small enough, we can construct 
\begin{eqnarray*}
u_{\chi_{t-\epsilon},b,t-\epsilon,N,\omega} & = & \chi_{t-\epsilon}(x')Q_{t-\epsilon}(x')e^{i\tau x\cdot\xi}e^{-\tau(x\cdot\omega-(t-\epsilon))A_{t-\epsilon}(x')}b+\gamma_{\chi_{t-\epsilon},b,t-\epsilon,N,\omega}\\
 &  & +r_{\chi_{t-\epsilon},b,t-\epsilon,N,\omega}
\end{eqnarray*}
to be the oscillating-decaying solution for $\nabla\cdot(A^{0}(x)\nabla\cdot)+k^{2}\cdot$
in $B_{t-\epsilon}(\omega)=B\cap\{x\cdot\omega>t-\epsilon\}$, where
$\chi_{t-\epsilon}(x')\in C_{0}^{\infty}(\mathbb{R}^{2})$ and $b\in\mathbb{C}$.
Note that in section 2, we have assumed the leading coefficient $A^{0}(x)\in B^{\infty}(\mathbb{R}^{3})$.
Similarly, we have the oscillating-decaying solution 
\[
u_{\chi_{t},b,t,N,\omega}(x,\tau)=\chi_{t}(x')Q_{t}e^{i\tau x\cdot\xi}e^{-\tau(x\cdot\omega-t)A_{t}(x')}b+\gamma_{\chi_{t},b,t,N,\omega}(x,\tau)+r_{\chi_{t},b,t,N,\omega}
\]
for $L_{A^{0}}$ in $B_{t}(\omega)$. In fact, for any $\tau$, $u_{\chi_{t-\epsilon},b,t-\epsilon,N,\omega}(x,\tau)\to u_{\chi_{t},b,t,N,\omega}(x,\tau)$
in an appropriate sense as $\epsilon\to0$. For details, we refer
readers to consult all the details and results in \cite{JN-ods},
and we list consequences in the following.
\[
\chi_{t-\epsilon}(x')Q_{t-\epsilon}(x')e^{i\tau x\cdot\xi}e^{-\tau(x\cdot\omega-(t-\epsilon))A_{t-\epsilon}(x')}b\to\chi_{t}(x')Q_{t}e^{i\tau x\cdot\xi}e^{-\tau(x\cdot\omega-t)A_{t}(x')}b
\]
in $H^{2}(B_{t}(\omega))$ as $\epsilon$ tends to 0,
\[
\gamma_{\chi_{t-\epsilon},b,t-\epsilon,N,\omega}\to\gamma_{\chi_{t},b,t,N,\omega}
\]
in $H^{2}(B_{t}(\omega))$ as $\epsilon$ tends to 0, and finally,
\[
r_{\chi_{t-\epsilon},b,t-\epsilon,N,\omega}\to r_{\chi_{t},b,t,N,\omega}
\]
in $H^{1}(B_{t}(\omega))$ as $\epsilon$ tends to 0.

Obviously, $B_{t-\epsilon}(\omega)$ is a convex set and $\overline{\Omega_{t}(\omega)}\subset B_{t-\epsilon}(\omega)$
for all $t\leq t_{0}$. By using the Runge approximation property,
we can see that there exists a sequence of functions $\tilde{u}_{\epsilon,j}$,
$j=1,2,\cdots$, such that 
\[
\tilde{u}_{\epsilon,j}\to u_{\chi_{t-\epsilon},b,t-\epsilon,N,\omega}\mbox{ in }H^{1}(B_{t}(\omega)),
\]
where $\tilde{u}_{\epsilon,j}\in H^{1}(\widetilde{\Omega})$ satisfy
$L_{A^{0}}\tilde{u}_{\epsilon,j}=0$ in $\widetilde{\Omega}$ for
all $\epsilon,j$. Define the indicator function $I(\tau,\chi_{t},b,t,\omega)$
by the formula:
\[
I(\tau,\chi_{t},b,t,\omega)=\lim_{\epsilon\to0}\lim_{j\to\infty}\int_{\partial}(\Lambda_{D}-\Lambda_{\emptyset})f_{\epsilon,j}\overline{f_{\epsilon,j}}d\sigma,
\]
where $f_{\epsilon,j}=\tilde{u}_{\epsilon,j}|_{\partial\Omega}$. 

Note that in \cite{JN-ods}, they assume that $D$ satisfying the
following condition: For each $\omega\in S^{2}$, there exist $c_{\omega}>0$,
$\epsilon_{\omega}>0$ and $p_{\omega}\in[0,1]$ such that 
\[
\dfrac{1}{c_{\omega}}s^{p_{\omega}}\leq\mu(\{x\in D|x\cdot\omega=t_{0}+s\})\leq c_{\omega}s^{p_{\omega}}\mbox{ for all }s\in(0,\epsilon_{\omega}),
\]
where $\mu$ is the surface measure, but we drop this condition in
the following theorem. Now the characterization of the convex hull
of $D$ is based on the following theorem:
\begin{thm}
(1) If $t<t_{0}$, then for any $\chi_{t}\in C_{0}^{\infty}(\mathbb{R}^{2})$
and $b\in\mathbb{C}$, we have 
\[
\limsup_{\tau\to\infty}|I(\tau,\chi_{t},b,t,\omega)|=0.
\]
.

(2) If $t=t_{0}$, then for any $\chi_{t_{0}}\in C_{0}^{\infty}(\mathbb{R}^{2})$
with $x_{0}'=(x_{0}\cdot\eta,x_{0}\cdot\zeta)$ being an interior
point of $\mbox{supp}(\chi_{t_{0}})$ and $0\neq b\in\mathbb{C}$,
we have 
\[
\liminf_{\tau\to\infty}|I(\tau,\chi_{t_{0}},b,t_{0},\omega)|>0.
\]
\end{thm}
\begin{proof}
(1) Note that we have a sequence of functions $\{\tilde{u}_{\epsilon,j}\}$
satisfies the equation $\nabla\cdot(A^{0}\nabla u)+k^{2}u=0\mbox{ in }\Omega,$
as in the beginning of the section 3, let $w_{\epsilon,j}=u-\tilde{u}_{\epsilon,j}$,
then $w_{\epsilon,j}$ satisfies the Dirichlet problem 
\[
\begin{cases}
\nabla\cdot(A(x)\nabla w_{\epsilon,j})+k^{2}w=-\nabla\cdot((\widetilde{A}\chi_{D}-A^{0}\chi_{D})\nabla\tilde{u}_{\epsilon,j}) & \mbox{ in }\Omega,\\
w_{\epsilon,j}=0 & \mbox{ on }\partial\Omega.
\end{cases}
\]
So we can apply (\ref{eq:3.14}) directly, which means 
\[
\int_{\partial\Omega}(\Lambda_{D}-\Lambda_{\emptyset})f_{\epsilon,j}\overline{f_{\epsilon,j}}d\sigma\leq C\|\tilde{u}_{\epsilon,j}\|_{W^{1,p}(D)}^{2}\leq C\|\tilde{u}_{\epsilon,j}\|_{H^{1}(D)}^{2},
\]
where the last inequality obtained by the H$\ddot{o}$lder's inequality.

By the Runge approximation property we have 
\[
\tilde{u}_{\epsilon,j}\to u_{\chi_{t-\epsilon},b,t-\epsilon,N,\omega}\mbox{ in }H^{1}(B_{t}(\omega))
\]
 as $j\to\infty$ and we know that the obstacle $D\subset B_{t}(\omega)$,
so we have 
\[
\|\tilde{u}_{\epsilon,j}-u_{\chi_{t-\epsilon},b,t-\epsilon,N,\omega}\|_{H^{1}(D)}\to0
\]
as $j\to\infty$ for all $\epsilon>0$. Moreover, we know that $u_{\chi_{t-\epsilon},b,t-\epsilon,N,\omega}\to u_{\chi_{t},b,t,N,\omega}$
as $\epsilon\to0$ in $H^{1}(B_{t}(\omega))$, which implies 
\[
\|\tilde{u}_{\epsilon,j}-u_{\chi_{t},b,t,N,\omega}\|_{H^{1}(D)}\to0
\]
as $\epsilon\to0$, $j\to\infty$. Now by the definition of $I(\tau,\chi_{t},b,t,\omega)$,
we have 
\[
I(\tau,\chi_{t},b,t,\omega)\leq C\|u_{\chi_{t},b,t,N,\omega}\|_{H^{1}(D)}^{2}.
\]
Now if $t<t_{0}$, we substitute $u_{\chi_{t},b,t,N,\omega}=w_{\chi_{t},b,t,N,\omega}+r_{\chi_{t},b,t,N,\omega}$
with $w_{\chi_{t},b,t,N,\omega}$ being described by (\ref{eq:2.3})
into 
\[
I(\tau,\chi_{t},b,t,\omega)\leq C(\int_{D}|u_{\chi_{t},b,t,N,\omega}|^{2}dx+\int_{D}|\nabla u_{\chi_{t},b,t,N,\omega}|^{2}dx)
\]
and use estimates (\ref{eq:2.4}), (\ref{eq:2.5}) to obtain that
\[
|I(\tau,\chi_{t},b,t,\omega)|\leq C\tau^{-2N-1}
\]
which finishes 
\[
\limsup_{\tau\to\infty}|I(\tau,\chi_{t},b,t,\omega)=0.
\]

For the second part, we use (\ref{eq:3.15}), which means that we
have 
\[
\int_{\partial\Omega}(\Lambda_{D}-\Lambda_{\emptyset})f_{\epsilon,j}\overline{f_{\epsilon,j}}d\sigma\geq c\int_{D}|\nabla\tilde{u}_{\epsilon,j}|^{2}dx-k^{2}\int_{\Omega}|\tilde{w}_{\epsilon,j}|^{2}dx-\int_{D}|\tilde{u}_{\epsilon,j}|^{2}dx.
\]
From (\ref{eq:4.1}) and the similar argument in the first part, it
is easy to get 
\begin{equation}
I(\tau,\chi_{t},b,t,\omega)\geq c\int_{D}|\nabla u_{\chi_{t},b,t,N,\omega}|^{2}dx-c\|u_{\chi_{t},b,t,N,\omega}\|_{W^{1,p}(D)}^{2},\label{eq:4.1}
\end{equation}
where $w_{\chi_{t},b,t,N,\omega}=u-u_{\chi_{t},b,t,N,\omega}$.
\end{proof}
For the remaining part, we need some extra estimates in the following
section.

\subsection{End of the proof of Theorem 4.1}

In view of the lower bound, we need to introduce the sets $D_{j,\delta}\subset D$,
$D_{\delta}\subset D$ in the following. Recall that $h_{D}(\omega)=\inf_{x\in D}x\cdot\omega$
and $t_{0}=h_{D}(\omega)=x_{0}\cdot\omega$ for some $x_{0}\in\partial D$.
$\forall\alpha\in\partial D\cap\{x\cdot\rho=h_{D}(\omega)\}:=K$,
define $B(\alpha,\delta)=\{x\in\mathbb{R}^{3};|x-\alpha|<\delta\}$
($\delta>0$). Note $K\subset\cup_{\alpha\in K}B(\alpha,\delta)$
and $K$ is compact, so there exists $\alpha_{1},\cdots,\alpha_{m}\in K$
such that $K\subset\cup_{j=1}^{m}B(\alpha_{j},\delta)$. Thus, we
define 
\[
D_{j,\delta}:=D\cap B(\alpha_{j},\delta)\mbox{ and }D_{\delta}:=\cup_{j=1}^{m}D_{j,\delta}.
\]
It is easy to see that 
\[
\int_{D\backslash D_{\delta}}e^{-p\tau(x\cdot\omega-t_{0})A_{t_{0}}(x')}bdx=O(e^{-pa\tau}),
\]
where $A_{t_{0}}(x')\in B^{\infty}(\mathbb{R}^{2})$ is bounded and
its real part strictly greater than 0. so $\exists a>0$ such that
$\mbox{Re}A_{t_{0}}(x')\geq a>0$. Let $\alpha_{j}\in K$, by rotation
and translation, we may assume $\alpha_{j}=0$ and the vector $\alpha_{j}-x_{0}=-x_{0}$
is parallel to $e_{3}=(0,0,1)$. Therefore, we consider the change
of coordinates near each $\alpha_{j}$ as follows:
\[
\begin{cases}
y'=x'\\
y_{3}=x\cdot\omega-t_{0},
\end{cases}
\]
where $x=(x_{1},x_{2},x_{3})=(x',x_{3})$ and $y=(y_{1},y_{2},y_{3})=(y',y_{3})$.
Denote the parametrization of $\partial D$ near $\alpha_{j}$ by
$l_{j}(y')$, then we have the following estimates.
\begin{lem}
For $q\leq2$, we have 
\begin{eqnarray}
\int_{D}|u_{\chi_{t_{0}},b,t_{0},N,\omega}|^{q}dx & \leq & c\tau^{-1}\sum_{j=1}^{m}\iint_{|y'|<\delta}e^{-aq\tau l_{j}(y')}dy'+O(\tau^{-1}e^{-qa\delta\tau})\nonumber \\
 &  & +O(e^{-qa\tau})+O(\tau^{-3})+O(\tau^{-2N-1}),\label{eq:4.2}
\end{eqnarray}

\begin{eqnarray}
\int_{D}|u_{\chi_{t_{0}},b,t_{0},N,\omega}|^{2}dx & \geq & C\tau^{-1}\sum_{j=1}^{m}\iint_{|y'|<\delta}e^{-2a\tau l_{j}(y')}dy'+O(\tau^{-1}e^{-2a\delta\tau})\nonumber \\
 &  & +O(\tau^{-3})+O(\tau^{-2N-1}),\label{eq:4.3}
\end{eqnarray}

\begin{eqnarray}
\int_{D}|\nabla u_{\chi_{t_{0}},b,t_{0},N,\omega}|^{q}dx & \leq & C\tau^{q-1}\sum_{j=1}^{m}\iint_{|y'|<\delta}e^{-qa\tau l_{j}(y')}dy'+O(\tau^{-1}e^{-aq\delta\tau})\nonumber \\
 &  & +O(e^{-qa\tau})+O(\tau^{-1})+O(\tau^{-2N-1}),\label{eq:4.4}
\end{eqnarray}
and 
\begin{eqnarray}
\int_{D}|\nabla u_{\chi_{t_{0}},b,t_{0},N,\omega}|^{2}dx & \geq & C\tau\sum_{j=1}^{m}\iint_{|y'|<\delta}e^{-2a\tau l_{j}(y')}dy'+O(\tau^{-1}e^{-2\delta a\tau})\nonumber \\
 &  & +O(\tau^{-1})+O(\tau^{-2N-1}).\label{eq:4.5}
\end{eqnarray}
\end{lem}
\begin{proof}
The proof follows from \cite{Sini}. We only prove (\ref{eq:4.2})
and (\ref{eq:4.3}) and the proof of (\ref{eq:4.4}) and (\ref{eq:4.5})
are similar arguments.

For (\ref{eq:4.2}): 
\begin{eqnarray*}
\int_{D}|u_{\chi_{t_{0}},b,t_{0},N,\omega}|^{q}dx & \leq & C\int_{D}e^{-qa\tau(x\cdot\omega-t_{0})}dx+C_{q}\int_{D}|\gamma_{\chi_{t_{0}},b,t_{0},N,\omega}|^{q}dx\\
 &  & +C_{q}\int_{D}|r_{\chi_{t_{0}},b,t_{0},N,\omega}|^{q}dx\\
 & \leq & C\int_{D_{\delta}}e^{-qa\tau(x\cdot\omega-t_{0})}dx+C\int_{D\backslash D_{\delta}}e^{-qa\tau(x\cdot\omega-t_{0})}dx\\
 &  & +C\int_{D}|\gamma_{\chi_{t_{0}},b,t_{0},N,\omega}|^{2}dx+C\int_{D}|r_{\chi_{t_{0}},b,t_{0},N,\omega}|^{2}dx\\
 & \leq & C\sum_{j=1}^{m}\iint_{|y'|<\delta}dy'\int_{l_{j}(y')}^{\delta}e^{-qa\tau y_{3}}dy_{3}+Ce^{-qa\tau}\\
 &  & +C\|\gamma_{\chi_{t_{0}},b,t_{0},N,\omega}\|_{L^{2}(D)}^{2}+C\|r_{\chi_{t_{0}},b,t_{0},N,\omega}\|_{H^{1}(D)}^{2}\\
 & \leq & C\tau^{-1}\sum_{j=1}^{m}\iint_{|y'|<\delta}e^{-aq\tau l_{j}(y')}dy'-\dfrac{C}{q}\tau^{-1}e^{-qa\delta\tau}\\
 &  & +Ce^{-qa\tau}+C\tau^{-3}+C\tau^{-2N-1}
\end{eqnarray*}
note that $D\subset\Omega_{t_{0}}(\omega)$, which proves (\ref{eq:4.1}).

For (\ref{eq:4.3}):
\begin{eqnarray*}
\int_{D}|u_{\chi_{t_{0}},b,t_{0},N.\omega}|^{2}dx & \geq & C\int_{D}e^{-2a\tau(x\cdot\omega-t_{0})}dx-C\|\gamma_{\chi_{t_{0}},b,t_{0},N,\omega}\|_{L^{2}(\Omega_{t_{0}}(\omega))}^{2}\\
 &  & -+C\|r_{\chi_{t_{0}},b,t_{0},N,\omega}\|_{H^{1}(\Omega_{t_{0}}(\omega))}^{2}\\
 & \geq & C\int_{D_{\delta}}e^{-2a\tau(x\cdot\omega-t_{0})}dx-C\tau^{-3}-C\tau^{-2N-1}\\
 & = & C\tau^{-1}\sum_{j=1}^{m}\iint_{|y'|<\delta}e^{-2a\tau l_{j}(y')}dy'-\dfrac{C}{2}\tau^{-1}e^{-2a\tau}\\
 &  & -C\tau^{-3}-C\tau^{-2N-1}.
\end{eqnarray*}

\end{proof}
Recall that we have (\ref{eq:4.1}), the lower bound of $I(\tau,\chi_{t_{0}},b,t_{0},\omega)$,
so we want to compare the order (in $\tau$) of $\|u_{\chi_{t_{0}},b,t_{0},N,\omega}\|_{L^{2}(D)}$,
$\|\nabla u_{\chi_{t_{0}},b,t_{0},N,\omega}\|_{L^{2}(D)}$, $\|u_{\chi_{t_{0}},b,t_{0},N,\omega}\|_{L^{p}(D)}$
and $\|\nabla u_{\chi_{t_{0}},b,t_{0},N,\omega}\|_{L^{p}(D)}$.
\begin{lem}
For $\max\{2-\epsilon,\dfrac{6}{5}\}<p\leq2$, we have the estimates
as follows:
\[
\dfrac{\|\nabla u_{\chi_{t_{0}},b,t_{0},N,\omega}\|_{L^{2}(D)}^{2}}{\|u_{\chi_{t_{0}},b,t_{0},N,\omega}\|_{L^{2}(D)}^{2}}\geq C\tau^{2},\mbox{ }\dfrac{\|u_{\chi_{t_{0}},b,t_{0},N,\omega}\|_{L^{p}(\Omega)}^{2}}{\|u_{\chi_{t_{0}},b,t_{0},N,\omega}\|_{L^{2}(D)}^{2}}\geq C\tau^{1-\frac{2}{p}}
\]
and 
\[
\dfrac{\|\nabla u_{\chi_{t_{0}},b,t_{0},N,\omega}\|_{L^{p}(D)}^{2}}{\|u_{\chi_{t_{0}},b,t_{0},N,\omega}\|_{L^{2}(D)}^{2}}\geq C\tau^{3-\frac{2}{p}}
\]
for $\tau\gg1$.\end{lem}
\begin{proof}
The idea of the proof comes from \cite{Sini}, but here we still need
to deal with the $\gamma_{\chi_{t_{0}},b,t_{0},N,\omega}$ and $r_{\chi_{t_{0}},b,t_{0},N,\omega}$
in $D\subset\Omega_{t_{0}}(\omega)$. Note that if $\partial D$ is
Lipschitz, in our parametrization $l_{j}(y')$, we have $l_{j}(y')\leq C|y'|$.
Hence, 
\begin{eqnarray*}
\sum_{j=1}^{m}\iint_{|y'|<\delta}e^{-2a\tau l_{j}(y')}dy' & \geq & C\sum_{j=1}^{m}\iint_{|y'|<\delta}e^{-2\tau|y'|}dy'\\
 & \geq & C\tau^{-1}\sum_{j=1}^{m}\iint_{|y'|<\tau\delta}e^{-2|y'|}dy'\\
 & = & O(\tau^{-1}).
\end{eqnarray*}
For simplicity, we define $u_{0}:=u_{\chi_{t_{0}},b,t_{0},N,\omega}$
in the following calculations. Using lemma 4.2, we obtain
\[
\dfrac{\int_{D}|\nabla u_{0}|^{2}dx}{\int_{D}|u_{0}|^{2}dx}
\]
 
\begin{eqnarray*}
 & \geq & C\dfrac{\tau\sum_{j=1}^{m}\iint_{|y'|<\delta}e^{-2a\tau l_{j}(y')}dy'+O(\tau^{-1}e^{-2a\delta\tau})+O(\tau^{-1})+O(\tau^{-2N-1})}{\tau^{-1}\sum_{j=1}^{m}\iint_{|y'|<\delta}e^{-2a\tau l_{j}(y')}dy'+O(\tau^{-1}e^{-2a\delta\tau})+O(\tau^{-3})+O(\tau^{-2N-1})}\\
 & \geq & C\tau^{2}\dfrac{1+\frac{O(\tau^{-2}e^{-2a\delta\tau})+O(\tau^{-2})+O(\tau^{-2N-2})}{\sum_{j=1}^{m}\iint_{|y'|<\delta}e^{-2a\tau l_{j}(y')}dy'}}{1+\frac{O(e^{-2a\delta\tau})+O(\tau^{-2})+O(\tau^{-2N})}{\sum_{j=1}^{m}\iint_{|y'|<\delta}e^{-2a\tau l_{j}(y')}dy'}}\\
 & = & O(\tau^{2})
\end{eqnarray*}
as $\tau\gg1$, where 
\[
\lim_{\tau\to\infty}\dfrac{O(\tau^{-2}e^{-2a\delta\tau})+O(\tau^{-2})+O(\tau^{-2N-2})}{\sum_{j=1}^{m}\iint_{|y'|<\delta}e^{-2a\tau l_{j}(y')}dy'}=0
\]
 and 
\[
\lim_{\tau\to\infty}\dfrac{O(e^{-2a\delta\tau})+O(\tau^{-2})+O(\tau^{-2N})}{\sum_{j=1}^{m}\iint_{|y'|<\delta}e^{-2a\tau l_{j}(y')}dy'}=0.
\]
Now, by using the H$\ddot{o}$lder's inequality with the exponent
$q=\dfrac{2}{p}\geq1$, we have 
\[
\sum_{j=1}^{m}\iint_{|y'|<\delta}e^{-pa\tau l_{j}(y')}dy'\leq C(\sum_{j=1}^{m}\iint_{|y'|<\delta}e^{-2a\tau l_{j}(y')}dy')^{\frac{p}{2}}.
\]
Hence we use lemma 4.2 again, we have 
\[
\dfrac{(\int_{D}|u_{0}|^{p}dx)^{\frac{2}{p}}}{\int_{D}|u_{0}|^{2}dx}
\]
\begin{eqnarray*}
 & \leq & C\dfrac{\tau^{-\frac{2}{p}}(\sum_{j=1}^{m}\iint_{|y'|<\delta}e^{-pa\tau l_{j}(y')}dy')^{\frac{2}{p}}+O(\tau^{-\frac{2}{p}}e^{-2a\delta\tau})+O(e^{-2a\tau})}{\tau^{-1}\sum_{j=1}^{m}\iint_{|y'|<\delta}e^{-2a\tau l_{j}(y')}dy'+O(\tau^{-1}e^{-2a\delta\tau})+O(\tau^{-3})+O(\tau^{-2N-1})}\\
 &  & +\dfrac{O(\tau^{-\frac{6}{p}})+O(\tau^{\frac{-4N-2}{p}})}{\tau^{-1}\sum_{j=1}^{m}\iint_{|y'|<\delta}e^{-2a\tau l_{j}(y')}dy'+O(\tau^{-1}e^{-2a\delta\tau})+O(\tau^{-3})+O(\tau^{-2N-1})}
\end{eqnarray*}
\begin{eqnarray*}
 & \leq & C\tau^{-\frac{2}{p}+1}\dfrac{\sum_{j=1}^{m}\iint_{|y'|<\delta}e^{-2a\tau l_{j}(y')}dy'+O(e^{-2c\delta\tau})+O(e^{-2a\tau}\tau^{\frac{2}{p}})}{\sum_{j=1}^{m}\iint_{|y'|<\delta}e^{-2a\tau l_{j}(y')}dy'+O(e^{-2a\delta\tau})+O(\tau^{-2})+O(\tau^{-2N})}\\
 &  & +\dfrac{O(\tau^{-\frac{4}{p}})+O(\tau^{\frac{-4N}{p}})}{\tau^{-1}\sum_{j=1}^{m}\iint_{|y'|<\delta}e^{-2a\tau l_{j}(y')}dy'+O(\tau^{-1}e^{-2a\delta\tau})+O(\tau^{-3})+O(\tau^{-2N-1})}
\end{eqnarray*}
\begin{eqnarray*}
 & = & \tau^{-\frac{2}{p}+1}\dfrac{1+\frac{O(e^{-2c\delta\tau})+O(e^{-2c\tau}\tau^{\frac{2}{p}})+O(\tau^{-\frac{4}{p}})+O(\tau^{\frac{-4N}{p}})}{\sum_{j=1}^{m}\iint_{|y'|<\delta}e^{-2a\tau l_{j}(y')}dy'}}{1+\frac{O(e^{-2a\delta\tau})+O(\tau^{-2})+O(\tau^{-2N})}{\sum_{j=1}^{m}\iint_{|y'|<\delta}e^{-2a\tau l_{j}(y')}dy'}}\\
 & = & O(\tau^{-\frac{2}{p}+1})
\end{eqnarray*}
as $\tau\gg1$ and 
\[
\dfrac{(\int_{D}|\nabla u_{0}|^{p}dx)^{\frac{2}{p}}}{\int_{D}|u_{0}|^{2}dx}
\]
\begin{eqnarray*}
 & \leq & C\dfrac{\tau^{(p-1)\frac{2}{p}}(\sum_{j=1}^{m}\iint_{|y'|<\delta}e^{-pa\tau l_{j}(y')}dy')^{\frac{2}{p}}+O(\tau^{-\frac{2}{p}}e^{-2a\delta\tau})+O(e^{-2a\tau})}{\tau^{-1}\sum_{j=1}^{m}\iint_{|y'|<\delta}e^{-2a\tau l_{j}(y')}dy'+O(\tau^{-1}e^{-2a\delta\tau})+O(\tau^{-3})+O(\tau^{-2N-1})}\\
 &  & +C\dfrac{O(\tau^{-\frac{2}{p}})+O(\tau^{\frac{-4N-2}{p}})}{\tau^{-1}\sum_{j=1}^{m}\iint_{|y'|<\delta}e^{-2a\tau l_{j}(y')}dy'+O(\tau^{-1}e^{-2a\delta\tau})+O(\tau^{-3})+O(\tau^{-2N-1})}
\end{eqnarray*}
\begin{eqnarray*}
 & \leq & C\tau^{3-\frac{2}{p}}\dfrac{\sum_{j=1}^{m}\iint_{|y'|<\delta}e^{-2a\tau l_{j}(y')}dy'+O(\tau^{-1}e^{-2a\delta\tau})+O(e^{-2a\tau}\tau^{\frac{2}{p}-1})}{\sum_{j=1}^{m}\iint_{|y'|<\delta}e^{-2a\tau l_{j}(y')}dy'+O(e^{-2a\delta\tau})+O(\tau^{-2})+O(\tau^{-2N})}\\
 &  & +C\dfrac{O(\tau^{-1})+O(\tau^{\frac{-4N}{p}-1})}{+O(\tau^{-\frac{2}{p}})+O(\tau^{\frac{-4N-2}{p}})}
\end{eqnarray*}
\begin{eqnarray*}
 & = & C\tau^{3-\frac{2}{p}}\dfrac{1+\frac{O(\tau^{-1}e^{-2a\delta\tau})+O(e^{-2a\tau}\tau^{\frac{2}{p}-1})+O(\tau^{-1})+O(\tau^{\frac{-4N}{p}-1})}{\sum_{j=1}^{m}\iint_{|y'|<\delta}e^{-2a\tau l_{j}(y')}dy'}}{1+\frac{O(e^{-2a\delta\tau})+O(\tau^{-2})+O(\tau^{-2N})}{\sum_{j=1}^{m}\iint_{|y'|<\delta}e^{-2a\tau l_{j}(y')}dy'}}\\
 & = & O(\tau^{3-\frac{2}{p}})
\end{eqnarray*}
as $\tau\gg1$. By (\ref{eq:4.1}) and above estimates, we have 
\begin{eqnarray*}
\dfrac{I(\tau,\chi_{t},b,t,\omega)}{\|u_{\chi_{t},b,t,N,\omega}\|_{L^{2}(D)}^{2}} & \geq & C\tau^{2}-C\tau^{1-\frac{2}{p}}-C\tau^{3-\frac{2}{p}}\\
 & \geq & C\tau^{2}
\end{eqnarray*}
for $\tau\gg1$. On the other hand, for $\|u_{\chi_{t},b,t,N,\omega}\|_{L^{2}(D)}$,
we have 
\begin{eqnarray*}
\int_{D}|u_{\chi_{t},b,t,N,\omega}|^{2}dx & \geq & C\tau^{-1}\sum_{j=1}^{m}\iint_{|y'|<\delta}e^{-2a\tau l_{j}(y')}dy'+O(\tau^{-1}e^{-qa\delta\tau})\\
 &  & +O(\tau^{-3})+O(\tau^{-2N-1})\\
 & \geq & C\tau^{-1}\sum_{j=1}^{m}\iint_{|y'|<\delta}e^{-2a\tau|y'|}dy'+O(\tau^{-1}e^{-qa\delta\tau})\\
 &  & +O(\tau^{-3})+O(\tau^{-2N-1})\\
 & \geq & C\tau^{-2}\sum_{j=1}^{m}\iint_{|y'|<\tau\delta}e^{-2a|y'|}dy'+O(\tau^{-1}e^{-qa\delta\tau})\\
 &  & +O(\tau^{-3})+O(\tau^{-2N-1})\\
 & = & O(\tau^{-2}).
\end{eqnarray*}
Therefore, we have 
\[
I(\tau,\chi_{t},b,t,\omega)\geq C\tau^{2}\|u_{\chi_{t},b,t,N,\omega}\|_{L^{2}(D)}^{2}\geq C>0
\]
for $\tau\gg1$.
\end{proof}
In view of theorem 4.1 and lemma 4.2, we can give an algorithm for
reconstructing the convex hull of an inclusion $D$ by the Dirichlet-to-Neumann
map $\Lambda_{D}$ as long as $A(x)$ and $D$ satisfy the described
conditions.\\
\textbf{Reconstruction algorithm}.
\begin{enumerate}
\item Give $\omega\in S^{2}$ and choose $\eta,\zeta,\xi\in S^{2}$ so that
$\{\eta,\zeta,\xi\}$ forms a basis of $\mathbb{R}^{3}$ and $\xi$
lies in the span of $\eta$ and $\zeta$;
\item Choose a starting $t$ such that $\Omega\subset\{x\cdot\omega\geq t\}$;
\item Choose a ball $B$ such that the center of $B$ lies on $\{x\cdot\omega=s\}$
for some $s<t$ and $\Omega\subset\overline{B_{t}(\omega)}$ and take
$0\neq b\in\mathbb{C}$;
\item Choose $\chi_{t}\in C_{0}^{\infty}(\mathbb{R}^{2})$ such that $\chi_{t}>0$
in $\Sigma_{t}(\omega)$ and $\chi_{t}=0$ on $\partial\Sigma_{t}(\omega)$;
\item Construct the oscillating-decaying solution $u_{\chi_{t-\epsilon},b,t-\epsilon,N,\omega}$
in $B_{t-\epsilon}(\omega)$ with $\chi_{t-\epsilon}=\chi_{t}$ and
the approximation sequence $\tilde{u}_{\epsilon,j}$ in $\widetilde{\Omega}$;
\item Compute the indicator function $I(\tau,\chi_{t},b,t,\omega)$ which
is determined by boundary measurements;
\item If $I(\tau,\chi_{t},b,t,\omega)\to0$ as $\tau\to\infty$, then choose
$t'>t$ and repeat (iv), (v), (vi);
\item If $I(\tau,\chi_{t},b,t,\omega)\nrightarrow0$ for some $\chi_{t'}$,
then $t'=t_{0}=h_{D}(\omega)$;
\item Varying $\omega\in S^{2}$ and repeat (i) to (viii), we can determine
the convex hull of $D$.\end{enumerate}


\begin{thebibliography}{10}
\bibitem{Lax} P. Lax, A stability theorem for solutions of abstract
differential equations and its application to the study of local behavior
of solutions of elliptic equations, Comm. Pure Appl. Math. 9 (1956)
747-766

\bibitem{JN-ods} G. Nakamura, G. Uhlmann, J.-N. Wang, Oscillating-decaying
solutions, Runge approximation property for the anisotropic elasticity
system and their applications to inverse problems, J. Math. Pures
Appl. 84 (2005) 21-54.

\bibitem{UhlWang-disCGO}G. Uhlmann and J.-N. Wang, Reconstructing
discontinuities using complex geometrical optics solutions, SIAM J.
Appl. Math., 68 (2008), pp. 1026-1044.

\bibitem{Majza} V. Maz'ja, Sobolev spaces, Springer-Verlag, Berling-Heidelberg-NewYork-Tokyo,
1985.

\bibitem{Meyer} N. Meyers. An $L^{p}$-estimate for the gradient
of solutions of second order elliptic divergence equations. Ann. Scuola
Norm. Sup. Pisa (3), 17 (1963), 189-206.

\bibitem{Sini} M. Sini and K. Yoshida, On the reconstruction of interfaces
using CGO solutions for the acoustic case, Preprint.

\bibitem{Nakamura} G. Nakamura, Applications of the oscillating-decaying
solutions to inverse problems, Preprint.

\bibitem{JN-acoustic} S. Nagayusu, G. Uhlmann, Jenn-Nan Wang, Reconstruction
of penetrable obstacles in acoustic scattering, SIAM J. Math. Anal.
Vol. 43, No.1, 189-211.

\bibitem{Ikehata} M. Ikehata, Reconstruction of the shape of the
inclusion by boundary measurements. Comm. Partial Differential Equations
23 (1998), 1231-1241.

\bibitem{Ike-enclos1} M. Ikehata, The enclosure method and its applications,
in: Analytic Extension Formulas and Their Applications (Fukuoka, 1999
/ Kyoto, 2000), in: Anal. Appl. Comput., vol. 9, Kluwer, Dordrecht,
2001, pp. 87-103.

\bibitem{Ike-enclosue} M. Ikehata, Enclosing a polygonal cavity in
a two-dimensional bounded domain from Cauchy data. Inverse Problems
15 (1999), 1231-1241.

\bibitem{Sylvester} J. Sylvester, G. Uhlmann, A global uniqueness
theorem for an inverse problems for an inverse boundary value problem,
Ann. of Math. (2) 125 (1987) 153-169.

\bibitem{Taylor} M. Taylor, Pseudodifferential Operators, Princeton
Univ. Press, Princeton, NJ, 1981. 

\bibitem{Malgrange} Existence et approximation des solutions des
$\acute{e}$quations aux d$\acute{e}$riv$\acute{e}$es partielles
et des $\acute{e}$quations de convolution, Ann. Inst. Fourier (Grenoble)
6 (1955-1956) 271-355. \end{thebibliography}
\end{document}